\newif\iflatexe
\latexefalse

\iflatexe
  \documentclass[12pt]{article}
  \usepackage{latexsym,a4}
\else
  \documentstyle[12pt,a4]{article}
\fi
\makeatletter
\@addtoreset{equation}{subsection}

\makeatother

\newcommand{\separate}{\medskip\noindent}

\def\newsection{ \separate
   \refstepcounter{subsection} 
   {\large\bf \thesubsection\kern.3em}
}
\def\mytheorem#1{
   \separate{\large\bf Theorem#1:\kern.3em}} 

\def\mylemma#1{
   \separate{\large\bf Lemma#1:\kern.3em}} 

\def\mycorollary#1{
   \separate{\large\bf Corollary#1:\kern.3em}} 

\def\myprop#1{
   \separate{\large\bf Proposition#1:\kern.3em}} 

\def\myremark#1{
   \separate{\large\bf Remark#1:\kern.3em}} 

\def\mydefinition#1{
   \separate{\large\bf Definition#1:\kern.3em}} 

\def\myconjecture#1{
   \separate{\large\bf Conjecture#1:\kern.3em}} 

\def\Proof{\separate\underline{Proof:}\kern1em}

\newcommand{\REELL}{{\setlength{\unitlength}{1em}
                     \begin{picture}(0.75,1)
                     \put(0,0){\line(0,1){0.69}}
                     \put(0,0){\sf R}
                     \end{picture}
                   }}

\newcommand{\FREELL}{\mbox{\tiny{\setlength{\unitlength}{1em}
                     \begin{picture}(0.6,0.5)
                     \put(0,0){\line(0,1){0.48}}
                     \put(0,0){\rm R}
                     \end{picture}
                   }}}
\newcommand{\Brr}{{\mathchoice{\REELL}{\REELL}{\!\FREELL}{\!\FREELL}}}

\def\QED{\hfill$\Box$}

\def\threespace{{\Brr^{\!\lower2pt\hbox{\mbox{\rm\scriptsize{3}}}}}}

\begin{document}

\renewcommand{\thefootnote}{\fnsymbol{footnote}}

\begin{center}
{\LARGE On a theorem by do~Carmo and Dajczer}

\vskip1cm
\begin{minipage}{6cm}
\begin{center}
G.~Haak\footnotemark[1]\\ Fachbereich Mathematik\\ TU Berlin\\ D-10623 Berlin
\end{center}
\end{minipage}
\vspace{0.5cm}

\end{center}

\footnotetext[1]{supported by Sonderforschungsbereich 288}

\begin{abstract}
We give a new proof of a theorem by M.P.~do Carmo and M.~Dajczer on
helicoidal surfaces of constant mean curvature.
\end{abstract}

\section{Introduction}
Let $G$ be a one-parameter group of proper Euclidean motions of $\Brr^3$
of the form
$$
g_t(x,y,z)=(x\cos t+y\sin t, -x\sin t+y\cos t, z+ht), t\in\Brr.
$$
I.e., $G$ is a group of helicoidal transformations with pitch
$h\in\Brr$. In the degenerate case $h=0$, $G$
becomes a group of pure rotations. Up to an affine change of
coordinates and reparametrization, all one-parameter
groups of Euclidean motions are either of this form or are groups of
pure translations.

In 1982 do~Carmo and Dajczer~\cite{doCDa:1} investigated surfaces of
constant mean curvature (CMC-surfaces) which are generated from a
plane curve by the action of a helicoidal group in the same way as a
rotational surface is generated by the action of a group of rotations.
They proved the following theorem:

\mytheorem{ 1} {\em
A complete immersed CMC-surface is helicoidal if
and only if it is in the associated family of a Delaunay surface.
}

\separate
They proved this result by introducing for each helicoidal CMC-immersion
the $2$-parameter family of helicoidal surfaces given by Bour's
Lemma~\cite{Bour:1} and evaluating the constant curvature condition
for the elements of these families.
This approach on one hand gives an explicit parametrization of
helicoidal CMC-immersions. On the other hand, it reaches its goal, the
proof of Theorem~1, in a fairly indirect way.

Since helicoidal surfaces still spawn interest~\cite{KoTa:1,Sasai:1}, 
we want to show in this note how Theorem~1 can be
obtained in a much simpler way using a more recent theorem of
Smyth~\cite{Sm:2} and some results from~\cite{DoHa:2}.
We state Smyth's theorem in the language of~\cite{DoHa:2}:

\mytheorem{ 2} {\em
Let $\Phi:M\rightarrow\Brr^3$, $M$ a Riemann
surface, be a complete conformally immersed
CMC-surface admitting a one-parameter group of self-isometries. Then
the simply connected cover of $M$ is the complex plane 
and the surface is either in the associated family of a
Delaunay surface or its metric is rotationally invariant.
}

\separate
Here a self-isometry is an automorphism of the simply connected cover
${\mathcal D}$ of $M$, which preserves the metric of the universal
covering immersion $\Psi:{\mathcal D}\rightarrow\Brr^3$ given by pulling
back the immersion $\Phi$ to ${\mathcal D}$. For details
see~\cite{DoHa:2}.  Those CMC-surfaces which have a rotationally
invariant metric are now commonly called Smyth surfaces.

We also introduce the notion of a space symmetry of a CMC-immersion
$\Phi:M\rightarrow\Brr^3$.  A space symmetry of $\Phi$ is a Euclidean
motion in $\Brr^3$ which preserves the image of $\Phi$ as a set.  The
relation between space symmetries and self-isometries was also studied
exhaustively in \cite{DoHa:2}.  By~\cite[Lemma~2.15]{DoHa:2} the group
of space symmetries of a Smyth surface is discrete.

\section{The proof of Theorem~1}
We can now give the proof of Theorem~1 right away:

\Proof
For a given CMC-immersion
$\Phi:M\rightarrow\Brr^3$ there exists (see
e.g.~\cite[Theorem~2.2]{DoHa:2}) a conformal structure on $M$ such
that $M$ becomes a Riemann surface and $\Phi$ becomes a conformal
CMC-immersion.
If $\Phi$ is also complete and in addition admits a one-parameter
group of helicoidal space symmetries, then
by~\cite[Prop.~2.12]{DoHa:2} and \cite[Corollary~2.6]{DoHa:2},
$\Phi$ admits also a one-parameter group of self-isometries. In particular it
satisfies the assumptions of Theorem~2 above. Since a group of
helicoidal Euclidean motions is never discrete, the surface cannot be
a Smyth surface. It therefore has to
be in the associated family of a Delaunay surface.

Conversely, by~\cite[Lemma~2.15]{DoHa:2} and~\cite[Prop.~3.4]{DoHa:2} 
each element of the associated family of a Delaunay
surface admits a one-parameter group of space symmetries. Since the
most general one-parameter group of Euclidean motions is a group
of helicoidal transformations (with possibly degenerate pitch), all
surfaces in the associated family of a Delaunay surface are helicoidal
or rotational.
\centerline{\hfill\QED}

It should also be noted that in the language of integrable systems
(the metric of a conformal CMC-immersion without umbilics satisfies the
integrable $\sinh$-Gordon equation), Theorem~1 also implies
that helicoidal CMC-surfaces are of finite type (see~\cite{PiSt:1} and
\cite{Bo:1}). Thus for helicoidal surfaces, 
apart from the parametrizations given
in~\cite{doCDa:1} and \cite{Sasai:1}, there is Bobenko's
parametrization in terms of theta functions~\cite{Bo:1}.

\end{document}